\documentclass[11pt]{amsart}
\usepackage{graphicx}
\usepackage[centertags]{amsmath}
\usepackage{amsfonts}
\usepackage{amssymb}
\usepackage{amsthm}
\usepackage{newlfont}
\newtheorem{theorem}{Theorem}[section]
\newtheorem{corollary}[theorem]{Corollary}
\newtheorem{lemma}[theorem]{Lemma}
\newtheorem{prop}[theorem]{Proposition}
\theoremstyle{definition}

\theoremstyle{remark}
\newtheorem{rem}[theorem]{Remark}
\newcommand{\tr}{\textrm{tr}}
\newcommand{\id}{\textrm{id}}

\title{Four-dimensional Riemannian  product manifolds with circulant structures}
\author{Iva Dokuzova}
\address{Iva Dokuzova\\Department of Algebra and Geometry\\
Plovdiv University Paisii Hilendarski\\ 24 Tzar Asen, 4000 Plovdiv, Bulgaria}

\email{dokuzova@uni-plovdiv.bg}

\begin{document}

\begin{abstract}
A 4-dimensional Riemannian manifold equipped with an additional tensor structure, whose fourth power is the identity, is considered.
This structure has a circulant matrix with respect to some basis, i.e. the structure is circulant, and it acts as an isometry with respect to the metric.
 The Riemannian product manifold associated with the considered manifold is studied.
 Conditions for the metric, which imply that the Riemannian product manifold belongs to each of the basic classes of Staikova-Gribachev's classification, are obtained.
 Examples of such manifolds are given.
\end{abstract}



\maketitle

\section{Introduction}\label{sec1}

The study of Riemannian manifolds $M$ with a metric $g$ and an almost product structure $P$ is initiated by K.~Yano in \cite{yano1}. The classification of the almost product manifolds $(M, g, P)$ with respect to the covariant derivative of $P$ is made by A.~M.~Naveira in \cite{nav}.
The manifolds $(M, g, P)$ with zero trace of $P$ are classified with respect to the covariant derivative of $P$ by M.~Staikova and K.~Gribachev in \cite{S-G}.
The basic classes of this classification are $\mathcal{W}_{1}, \mathcal{W}_{2}$ and $\mathcal{W}_{3}$. Their intersection is the class $\mathcal{W}_{0}$ of Riemannian $P$-manifolds. The class of the Riemannian product manifolds is $\mathcal{W}_1\oplus\mathcal{W}_2$. It is formed by manifolds with an integrable structure $P$. The class $\mathcal{W}_1$ consists of the conformal Riemannian $P$-manifolds. Some of the recent studies of the Riemannian almost product manifolds are made in \cite{druta},  \cite{gribacheva}, \cite{gribacheva3}, \cite{griba-mek} and \cite{pusic}.

Problems of differential geometry of a $4$-dimensional Riemannian manifold $M$ with a tensor structure $Q$ of type $(1,1)$, which satisfies $Q^{4}=\id$, $Q^{2}\neq\pm\id$, are considered in \cite{dokuzova}, \cite{dokuzova2}, \cite{dok-raz-man}, \cite{G-R} and \cite{raz}. The matrix of $Q$ in some basis is circulant and $Q$ is compatible with the metric $g$, so that an isometry is induced in any tangent space on $M$. Such a manifold $(M, g, Q)$ is associated with a Riemannian almost product manifold $(M, g, P)$, where $P=Q^{2}$, $\tr P=0$.

 In the present work we continue studying the manifold $(M, g, Q)$ and the associated manifold $(M, g, P)$. Our purpose is to determine their position, among the well-known manifolds, by using the classifications in \cite{nav} and \cite{S-G}. In Sect.~\ref{2}, we recall some necessary facts about $(M, g, Q)$ and about $(M, g, P)$. In Sect.~\ref{3}, we obtain the components of the fundamental tensor $F$, determined by the metric $g$ and by the covariant derivative of $P$. We establish that  $(M, g, P)$ is a Riemannian product manifold, i.e. $(M, g, P)$ belongs to $\mathcal{W}_1\oplus\mathcal{W}_2$. We find necessary and sufficient conditions that $(M, g, P)$ belongs to each of the classes $\mathcal{W}_0$, $\mathcal{W}_1$ and $\mathcal{W}_2$. In Sect.~\ref{4}, we give some examples of the considered Riemannian product manifolds.
 \section{Preliminaries}\label{2}
  Let $M$ be a $4$-dimensional Riemannian manifold equipped with a tensor structure $Q$ in the tangent space $T_{p}M$ at
an arbitrary point $p$ on $M$. The structure $Q$ has a circulant matrix, with respect to some basis, as follows:
\begin{equation}\label{Q}
    (Q_{j}^{k})=\begin{pmatrix}
      0 & 1 & 0 & 0\\
      0 & 0 & 1 & 0\\
      0 & 0 & 0 & 1\\
      1 & 0 & 0 & 0\\
    \end{pmatrix}.
\end{equation}
Then $Q$ has the properties
\begin{equation*}
    Q^{4}=\id,\quad Q^{2}\neq\pm\id.
\end{equation*}
Let the metric $g$ and the structure $Q$ satisfy
\begin{equation*}
 g(Qx, Qy)=g(x,y),\quad x, y\in \mathfrak{X}(M).
\end{equation*}
The above condition and \eqref{Q} imply that the matrix of $g$ has the form:
\begin{equation}\label{metricg}
    (g_{ij})=\begin{pmatrix}
      A & B & C & B \\
      B & A & B & C\\
      C & B & A & B\\
      B & C & B & A\\
    \end{pmatrix}.
\end{equation}
Here $A$, $B$ and $C$ are smooth functions of an arbitrary point $p(x^{1}, x^{2},
x^{3}, x^{4})$ on $M$.
It is supposed that $A > C > B >0$ in order $g$ to be positive definite. The manifold $(M, g, Q)$ is introduced in \cite{raz}.

Anywhere in this work, $x, y, z$ will stand for arbitrary elements of the algebra on the smooth vector fields on $M$ or vectors in the tangent space $T_{p}M$. The Einstein summation convention is used, the range of the summation indices being always $\{1, 2, 3, 4\}$.

 In \cite{dokuzova}, it is noted that the manifold $(M, g, P)$, where $P = Q^{2}$, is a Riemannian manifold with an almost product structure $P$, because $P^{2} = \id $, $P\neq\pm \id$ and $g(Px, Py) = g(x, y)$. Moreover $\tr P = 0$. For such manifolds Staikova-Gribachev's classification is valid \cite{S-G}.
This classification was made with respect to the tensor $F$ of type $(0,3)$ and the associated 1-form $\theta$, which are defined by
\begin{equation}\label{F}
  F(x,y,z)=g((\nabla_{x}P)y,z),\quad \theta(x)=g^{ij}F(e_{i},e_{j},x).
\end{equation}
Here $\nabla$ is the Levi-Civita connection of $g$, and $g^{ij}$ are the components of the inverse matrix of $(g_{ij})$ with respect to an arbitrary basis $\{e_{i}\}$.
The tensor $F$ has the following properties
\begin{equation}\label{prop-F}
  F(x,z,y)=F(x,y,z),\quad F(x, Py, Pz)=-F(x, y, z).
\end{equation}

The manifolds $(M, g, P)$ with integrable structure $P$ are called Riemannian product manifolds and they form the class  $\mathcal{W}_1\oplus\mathcal{W}_2$. The characteristic conditions for the classes $\mathcal{W}_{0}$, $\mathcal{W}_{1}, \mathcal{W}_{2}$ and  $\mathcal{W}_1\oplus\mathcal{W}_2$ are the following:
\begin{enumerate}
  \item[i)] $\mathcal{W}_{0}$:
  \begin{equation}\label{c0}
 \quad F(x, y, z)=0,
\end{equation}
  \item[ii)] $\mathcal{W}_{1}$:
  \begin{equation}\label{c1}
\begin{split}
 F(x,y,z)=\frac{1}{4}\big[&(g(x,y)\theta(z)+g(x,z)\theta(y)\\&-g(x,Py)\theta(Pz)-g(x,Pz)\theta(Py)\big],
\end{split}
\end{equation}
  \item[iii)] $\mathcal{W}_{2}$:
  \begin{equation}\label{c2}
\begin{split}
 F(x,y,Pz)+F(y,z,Px)+F(z,x,Py)=0,\quad \theta(z)=0,
\end{split}
\end{equation}
  \item[iv)] $\mathcal{W}_{1}\oplus\mathcal{W}_{2}$:
  \begin{equation}\label{c3}
\begin{split}
 F(x,y,Pz)+F(y,z,Px)+F(z,x,Py)=0.
\end{split}
\end{equation}
\end{enumerate}

In the next section, we obtain conditions that $(M, g, P)$ belongs to each of the classes above and to some of their subclasses. Thus
the following statements are useful for the completeness of our research.
\begin{theorem}\cite{G-R} If the structure $Q$ of $(M, g, Q)$ satisfies $\nabla Q=0$, then $\nabla P=0$, i.e. $(M, g, P)$ belongs to the class $\mathcal{W}_{0}$.
\end{theorem}
\begin{theorem}\cite{raz}
The structure $Q$ of $(M, g, Q)$ satisfies $\nabla Q=0$ if and only if the following equalities are valid:
\begin{equation}\label{tQ}
\begin{array}{llll}
A_{3}=C_{1}, & A_{1}=C_{3}, &  B_{3}=B_{1}, &  2B_{1}=C_{4}+C_{2},\\
B_{4}=B_{2}, & A_{2}=C_{4}, & A_{4}=C_{2}, & 2B_{2}=C_{1}+C_{3},
\end{array}
\end{equation}
where $A_{i}=\frac{\partial A}{\partial x^{i}},$ $B_{i}=\frac{\partial B}{\partial x^{i}}$,  $C_{i}=\frac{\partial C}{\partial x^{i}}$.
\end{theorem}

\section{The fundamental tensor $F$ on $(M, g, P)$}\label{3}
The components of Nijenhuis tensor $N$ of the almost product structure $P$ are determined by the equalities
\begin{equation*}
N_{ij}^{k}=P_{i}^{a}(\partial_{a}P_{j}^{k}-\partial_{j}P_{a}^{k})-P_{j}^{a}(\partial_{a}P_{i}^{k}-\partial_{i}P_{a}^{k}).
\end{equation*}
It is known from \cite{S-G} that the vanishing of the Nijenhuis tensor $N$ is equivalent to the condition \eqref{c3}.

Using \eqref{Q}, we get that the components of the almost product structure $P=Q^{2}$ on $(M, g, P)$ are given by the matrix
\begin{equation}\label{P}
    (P_{j}^{k})=\begin{pmatrix}
      0 & 0 & 1 & 0\\
      0 & 0 & 0 & 1\\
      1 & 0 & 0 & 1\\
      0 & 1 & 0 & 0\\
    \end{pmatrix}.
\end{equation}
Evidently the Nijenhuis tensor $N$ of $P$ vanishes, so we have the following
\begin{lemma}\label{LW12}
The manifold $(M, g, P)$ belongs to the class $\mathcal{W}_{1}\oplus\mathcal{W}_{2}$.
\end{lemma}
Further, we will consider each of the cases when the fundamental tensor $F$ on $(M, g, P)$ satisfies the identities \eqref{c0}, \eqref{c1} or \eqref{c2}, since $F$ satisfies \eqref{c3}. First we calculate the components of $F$.
\begin{lemma}\label{tF}
 The nonzero components $F_{ijk}=F(e_{i}, e_{j}, e_{k})$ of the fundamental tensor $F$ on the manifold $(M, g, P)$ are given by
\begin{equation}\label{nablaf}
\begin{array}{ll}
F_{111}=-F_{133}=A_{3}-C_{1}, & 2F_{112}=-2F_{134}=A_{4}-B_{1}-C_{2}+B_{3},\\
F_{122}=-F_{144}=B_{4}-B_{2}, & 2F_{114}=-2F_{123}=A_{2}-B_{1}-C_{4}+B_{3},\\
F_{211}=-F_{233}=B_{3}-B_{1}, & 2F_{212}=-2F_{234}=A_{3}-B_{2}-C_{1}+B_{4},\\
F_{222}=-F_{244}=A_{4}-C_{2}, & 2F_{223}=-2F_{214}=A_{1}+B_{4}-C_{3}-B_{2},\\
F_{311}=-F_{333}=C_{3}-A_{1}, & 2F_{334}=-2F_{312}=A_{2}+B_{1}-C_{4}-B_{3},\\
F_{322}=-F_{344}=B_{4}-B_{2}, & 2F_{323}=-2F_{314}=A_{4}+B_{1}-C_{2}-B_{3},\\
F_{411}=-F_{433}=B_{3}-B_{1}, & 2F_{434}=-2F_{412}=A_{1}+B_{2}-C_{3}-B_{4},\\
F_{422}=-F_{444}=C_{4}-A_{2}, & 2F_{414}=-2F_{423}=A_{3}-B_{4}-C_{1}+B_{2}.
\end{array}
\end{equation}
\end{lemma}
\proof
The inverse matrix of $(g_{ij})$ has the form:
\begin{equation}\label{g-obr}
    (g^{ik})=\frac{1}{D}\begin{pmatrix}
            \overline{A}& \overline{B}& \overline{C} &\overline{B}\\
            \overline{B}&\overline{A}&\overline{B} &\overline{C}\\
            \overline{C}&\overline{B}&\overline{A}&\overline{B}\\
            \overline{B} & \overline{C} &\overline{B} &\overline{A}\\
            \end{pmatrix},
  \end{equation}
where
\begin{equation}\label{coef-obr}
\begin{array}{ll}
\overline{A}=A(A+C)-2B^{2}, & \overline{B}=B(C-A),\\ \overline{C}=2B^{2}-C(A+C),  & D=(A-C)\big[(A+C)^{2}-4B^{2}\big].
\end{array}
\end{equation}

Let $\Gamma_{ij}^{s}$ be the Christoffel symbols of $\nabla$. The next formula is well known:
\begin{equation*}
2\Gamma_{ij}^{k}=g^{ak}(\partial_{i}g_{aj}+\partial_{j}g_{ai}-\partial_{a}g_{ij}).
\end{equation*}
Then, using \eqref{metricg}, \eqref{g-obr} and \eqref{coef-obr}, we calculate these coefficients. They are given below:
\begin{align}\label{gama}\nonumber
\Gamma_{ii}^{i}=\frac{1}{2D}\big[&B(C-A)(4B_{i}-A_{j}-A_{s})+2B^{2}(2C_{i}-A_{i}-A_{k})\\\nonumber &+A(A+C)A_{i}-C(A+C)(2C_{i}-A_{k})\big],\\\nonumber
\Gamma_{ii}^{j}=\frac{1}{2D}\big[&B(C-A)(A_{i}+2C_{i}-A_{k})+2B^{2}(A_{j}-A_{s})\\\nonumber &+A(A+C)(2B_{i}-A_{j})-C(A+C)(2B_{i}-A_{s})\big],\\\nonumber
\Gamma_{ii}^{k}=\frac{1}{2D}\big[&B(C-A)(4B_{i}-A_{j}-A_{s})+2B^{2}(A_{i}+A_{k}-2C_{i})\\\nonumber &+A(A+C)(2C_{i}-A_{k})-C(A+C)A_{i}\big],\\
\Gamma_{ij}^{i}=\frac{1}{2D}\big[&B(C-A)(A_{i}+C_{i}+B_{j}-B_{s})+2B^{2}(B_{i}+C_{j}-B_{k}-A_{j})\\\nonumber &+A(A+C)A_{j}-C(A+C)(B_{i}+C_{j}-B_{k})\big],\\\nonumber
\Gamma_{ij}^{k}=\frac{1}{2D}\big[&B(C-A)(A_{i}+C_{i}+B_{j}-B_{s})+2B^{2}(B_{k}-C_{j}-B_{i}+A_{j})\\\nonumber &+A(A+C)(B_{i}+C_{j}-B_{k})-C(A+C)A_{j}\big],\\\nonumber
\Gamma_{ik}^{j}=\frac{1}{2D}\big[&B(C-A)(A_{i}+A_{k})+2B^{2}(C_{j}-C_{s})\\\nonumber &+A(A+C)(B_{i}+B_{k}-C_{j})-C(A+C)(B_{i}+B_{k}-C_{s})\big],\\\nonumber
\Gamma_{ik}^{i}=\frac{1}{2D}\big[&B(C-A)(2B_{i}+2B_{k}-C_{j}-C_{s})+2B^{2}(A_{i}-A_{k})\\\nonumber &+A(A+C)A_{k}-C(A+C)A_{i})\big].
\end{align}
In the equalities \eqref{gama} it is assumed that $i\neq j\neq k\neq s$ and the numbers in the pair $(i, k)$ (resp. in $(j, s)$) are simultaneously even or odd.

The  matrix  of the associated metric $\tilde{g}$, determined by $\tilde{g}(x,y)=g(x, Py)$, is
of the type:
\begin{equation}\label{tilde-g}
(\tilde{g}_{ij})=\begin{pmatrix}
      C & B & A & B\\
      B & C & B & A\\
      A & B & C & B\\
      B & A & B & C\\
    \end{pmatrix}.
\end{equation}
Due to \eqref{F} the components of $F$ are $F_{ijk}=\nabla_{i}\tilde{g}_{jk}$.
The well-known are the following identities for a Riemannian metric:
\begin{equation}\label{defF}
\nabla_{i}\tilde{g}_{jk}=\partial_{i}\tilde{g}_{jk}-\Gamma_{ij}^{a}\tilde{g}_{ak}-\Gamma_{ik}^{a}\tilde{g}_{aj}.
\end{equation}
Applying \eqref{gama} and \eqref{tilde-g} into \eqref{defF}, and bearing in mind \eqref{F} and \eqref{prop-F}, we find the nonzero components of $F$, given in \eqref{nablaf}.
\endproof

Immediately, we have the following
\begin{corollary}\label{t-theta}
 The components $\theta_{k}=g^{ij}F(e_{i}, e_{j}, e_{k})$ of the 1-form $\theta$ on the manifold $(M, g, P)$ are expressed by the equalities
 \begin{equation}\label{theta}
\begin{split}
\theta_{1}&=\frac{1}{D}\big[\overline{C}(2C_{3}-2A_{1})+\overline{B}(4B_{3}-4B_{1})+\overline{A}(2A_{3}-2C_{1})\big],\\
\theta_{2}&=\frac{1}{D}\big[\overline{A}(2A_{4}-2C_{2})+\overline{B}(4B_{4}-4B_{2})+\overline{C}(2C_{4}-2A_{2})\big],\\
\theta_{3}&=\frac{1}{D}\big[\overline{C}(2C_{1}-2A_{3})+\overline{B}(4B_{1}-4B_{3})+\overline{A}(2A_{1}-2C_{3})\big],\\
\theta_{4}&=\frac{1}{D}\big[\overline{A}(2A_{2}-2C_{4})+\overline{B}(4B_{2}-4B_{4})+\overline{C}(2C_{2}-2A_{4})\big].
\end{split}
\end{equation}
 \end{corollary}
 \proof The equalities \eqref{theta} follow by direct computations from \eqref{nablaf}, \eqref{g-obr} and \eqref{coef-obr}.
 \endproof

Having in mind Lemma~\ref{LW12}, Lemma~\ref{tF} and Corollary~\ref{t-theta} we obtain the next statements.
\begin{theorem}\label{tw0}
The manifold $(M, g, P)$ belongs to the class $\mathcal{W}_0$ if and only if the following equalities are valid:
\begin{equation}\label{w0}
\begin{array}{lll}
A_{3}=C_{1}, & A_{1}=C_{3}, &  B_{3}=B_{1},\\
B_{4}=B_{2}, & A_{2}=C_{4}, & A_{4}=C_{2}.
\end{array}
\end{equation}
\end{theorem}
\proof
Due to \eqref{nablaf} we get that \eqref{c0} is satisfied if and only if \eqref{w0} holds true.
\endproof

\begin{theorem}\label{tw1}
The manifold $(M, g, P)$ belongs to the class $\mathcal{W}_1$ if and only if the following equalities are valid:
\begin{equation}\label{w1}
\begin{split}
(A+C)(B_{4}-B_{2})=B(A_{4}-C_{2}+C_{4}-A_{2}),\\ (A+C)(B_{3}-B_{1})=B(A_{3}-C_{1}+C_{3}-A_{1}).
\end{split}
\end{equation}
\end{theorem}
\proof
Using \eqref{metricg},  \eqref{P}, \eqref{nablaf},  \eqref{tilde-g}, \eqref{theta} and \eqref{w1} we obtain
 \begin{equation}\label{usl-w1}
 \begin{split}
  F_{kij}=\frac{1}{4}\big(g_{kj}\theta_{i}+g_{ki}\theta_{j}-\tilde{g}_{kj}\tilde{\theta}_{i}-\tilde{g}_{ki}\tilde{\theta}_{j}\big),\quad \tilde{\theta}_{i}=P_{i}^{a}\theta_{a},\\
  \end{split}
\end{equation}
which is equivalent to \eqref{c1}.

Vice versa, if \eqref{usl-w1} holds true, then \eqref{metricg}, \eqref{P}, \eqref{nablaf}, \eqref{tilde-g} and \eqref{theta} imply \eqref{w1}.
\endproof

In \cite{S-G}, it is proved that $\mathcal{W}_{1}=\mathcal{\overline{W}}_{3}\oplus\mathcal{\overline{W}}_{6}$, where $\mathcal{\overline{W}}_{3}$ and $\mathcal{\overline{W}}_{6}$ are two basic classes of Naveira's classification. These classes have the following characteristic conditions (\cite{gribacheva2}, \cite{griba-mek2}):
\begin{equation}\label{sub3c1}
 \begin{split}
  \mathcal{\overline{W}}_{3}:\quad F(x,y,z)=\frac{1}{4}\Big[&\big(g(x,y)+g(x,Py)\big)\theta(z)\\&+\big(g(x,z)+g(x,Pz)\big)\theta(y)\Big],\ \theta(Pz)=-\theta(z).
  \end{split}
\end{equation}
\begin{equation}\label{sub6c1}
 \begin{split}
  \mathcal{\overline{W}}_{6}:\quad F(x,y,z)=\frac{1}{4}\Big[&\big(g(x,y)-g(x,Py)\big)\theta(z)\\&+\big(g(x,z)-g(x,Pz)\big)\theta(y)\Big],\quad \theta(Pz)=\theta(z).
  \end{split}
\end{equation}
\begin{corollary}\label{tsubw3}
The manifold $(M, g, P)$ belongs to the class $\mathcal{\overline{W}}_{3}$ if and only if the following equalities are valid:
\begin{equation}\label{subw3}
\begin{array}{ll}
A_{4}-C_{2}=C_{4}-A_{2},& (A+C)(B_{4}-B_{2})=2B(A_{4}-C_{2}),\\
A_{3}-C_{1}=C_{3}-A_{1},& (A+C)(B_{3}-B_{1})=2B(A_{3}-C_{1}).
\end{array}
\end{equation}
\end{corollary}
\proof
The local form of \eqref{sub3c1} is
 \begin{equation}\label{usl-w3}
 \begin{split}
  F_{kij}=\frac{1}{4}\Big[\big(g_{ki}+\tilde{g}_{ki}\big)\theta_{j}+\big(g_{kj}+\tilde{g}_{kj}\big)\theta_{i}\Big],\quad \tilde{\theta}_{i}=-\theta_{i}.
  \end{split}
\end{equation}
Taking into account \eqref{theta} and \eqref{w1}, we get that \eqref{usl-w3} is satisfied if and only if \eqref{subw3} holds true.
\endproof
\begin{corollary}\label{tsubw6}
The manifold $(M, g, P)$ belongs to the class $\mathcal{\overline{W}}_{6}$ if and only if the following equalities are valid:
\begin{equation}\label{subw6}
\begin{array}{ll}
A_{4}-C_{2}=A_{2}-C_{4}, & B_{4}=B_{2}, \\
A_{3}-C_{1}=A_{1}-C_{3},& B_{3}=B_{1}.
\end{array}
\end{equation}
\end{corollary}
\proof
The local form of \eqref{sub6c1} is
 \begin{equation}\label{usl-w6}
 \begin{split}
  F_{kij}=\frac{1}{4}\Big[\big(g_{ki}-\tilde{g}_{ki}\big)\theta_{j}+\big(g_{kj}-\tilde{g}_{kj}\big)\theta_{i}\Big],\quad \tilde{\theta}_{i}=\theta_{i}.
  \end{split}
\end{equation}
Having in mind \eqref{theta} and \eqref{w1}, we get that \eqref{usl-w6} is satisfied if and only if \eqref{subw6} holds true.
\endproof
\begin{theorem}\label{tw2}
The manifold $(M, g, P)$ belongs to the class $\mathcal{W}_2$ if and only if the following equalities are valid:
\begin{equation}\label{w2}
\begin{split}
(A+C)(C_{3}-A_{1})=2B(B_{3}-B_{1}),\qquad C_{3}-A_{1}=A_{3}-C_{1},\\
(A+C)(C_{4}-A_{2})=2B(B_{4}-B_{2}),\qquad C_{4}-A_{2}=A_{4}-C_{2}.
 \end{split}
\end{equation}
\end{theorem}
\proof
From \eqref{c2} and  \eqref{theta} we have
\begin{equation}\label{23}
\begin{split}
\overline{C}(2C_{3}-2A_{1})+\overline{B}(4B_{3}-4B_{1})+\overline{A}(2A_{3}-2C_{1})=0,\\
\overline{A}(2A_{4}-2C_{2})+\overline{B}(4B_{4}-4B_{2})+\overline{C}(2C_{4}-2A_{2})=0,\\
\overline{C}(2C_{1}-2A_{3})+\overline{B}(4B_{1}-4B_{3})+\overline{A}(2A_{1}-2C_{3})=0,\\
\overline{A}(2A_{2}-2C_{4})+\overline{B}(4B_{2}-4B_{4})+\overline{C}(2C_{2}-2A_{4})=0.
\end{split}
\end{equation}
The equalities \eqref{coef-obr} and \eqref{23} imply \eqref{w2}.

Vice versa. We apply \eqref{w2} into \eqref{theta} and we get that \eqref{c2} holds true.
\endproof

\section{Examples of manifolds $(M, g, P)$}\label{4}
In this section we give a solution of each system of differential equations \eqref{w0}, \eqref{w1}, \eqref{subw3}, \eqref{subw6} and \eqref{w2},
in order to get examples of $(M, g, P)$ of the classes considered.
\subsection{An example in $\mathcal{W}_{0}$}
Let $(M, g, P)$ be a manifold with
\begin{equation}\label{prw0}
\begin{split}
    A&=(x^{1})^{2}+(x^{2})^{2}+(x^{3})^{2}+(x^{4})^{2},\\ B&=x^{1}+x^{2}+x^{3}+x^{4} ,\quad C=2x^{1}x^{3}+2x^{2}x^{4},
    \end{split}
\end{equation}
where $x^{i}>1.$

Evidently $A>C>B>0$ are valid. We check directly that the functions \eqref{prw0} and their derivatives satisfy the equalities \eqref{w0}.

Thus we have the following
\begin{prop}\label{kt2}
 The manifold $(M, g, P)$ with \eqref{prw0} belongs to $\mathcal{W}_{0}$.
\end{prop}

 \begin{rem} We note that the functions \eqref{prw0} do not satisfy \eqref{tQ}. Then we have $\nabla Q\neq 0$ for $(M, g, P)$, where $P=Q^{2}$.
 An example of a manifold $(M, g, Q)$ with $\nabla Q=0$ is given in \cite{raz}.
 \end{rem}

\subsection{An example in $\mathcal{\overline{W}}_{3}$}
Let $(M, g, P)$ be a manifold with
\begin{equation}\label{prw3}
\begin{split}
    A&=a(x^{1}+x^{2}-x^{3}-x^{4}),\quad B=b(x^{1}+x^{2}-x^{3}-x^{4}), \\ C&=c(x^{1}+x^{2}-x^{3}-x^{4}),
    \end{split}
\end{equation}
where $a, c, b \in \mathbb{R},\ a>c>b>0,\ x^{1}+x^{2}-x^{3}-x^{4}>0.$

The inequalities $A>C>B>0$ hold true.
The functions \eqref{prw3} and their derivatives satisfy the equalities \eqref{subw3} and do not satisfy the conditions \eqref{w0}.

Therefore, we establish the following
\begin{prop}
 The manifold $(M, g, P)$ with \eqref{prw3} belongs to $\mathcal{\overline{W}}_{3}$ but does not belong to $\mathcal{W}_{0}$.
\end{prop}

\subsection{An example in $\mathcal{\overline{W}}_{6}$}
Let $(M, g, P)$ be a manifold with
\begin{equation}\label{prw6}
\begin{split}
    A&=a(x^{1}+x^{2}+x^{3}+x^{4}),\quad B=b(x^{1}+x^{2}+x^{3}+x^{4}), \\ C&=c(x^{1}+x^{2}+x^{3}+x^{4}),
    \end{split}
\end{equation}
where $a, c, b \in \mathbb{R},\ a>c>b>0,\ x^{1}+x^{2}+x^{3}+x^{4}>0$.

The inequalities $A>C>B>0$ are satisfied.
The functions \eqref{prw6} and their derivatives satisfy the equalities \eqref{subw6} and do not satisfy the conditions \eqref{w0}.

Immediately, we state the following
\begin{prop}
 The manifold $(M, g, P)$ with \eqref{prw6} belongs to $\mathcal{\overline{W}}_{6}$ but does not belong to $\mathcal{W}_{0}$.
\end{prop}

\subsection{An example in $\mathcal{W}_{1}$}
Let $(M, g, P)$ be a manifold with
\begin{equation}\label{prw1}
\begin{split}
A=&a\exp(x^{1}-x^{2}),\ C=c\exp(x^{4}-x^{3}),\\
B=&a\exp(x^{1}-x^{2})-c\exp(x^{4}-x^{3}),
    \end{split}
\end{equation}
where $a, c \in \mathbb{R}^{+}$, $\ln\dfrac{c}{a}<x^{1}-x^{2}+x^{3}-x^{4}<\ln\dfrac{2c}{a}.$

Then $A>C>B>0$ are valid. The functions \eqref{prw1} satisfy the equalities \eqref{w1} but do not satisfy any of the conditions \eqref{subw3} and \eqref{subw6}.

Therefore, we establish the following
\begin{prop}
 The manifold $(M, g, P)$ with \eqref{prw1} belongs to $\mathcal{W}_{1}$ but does not belong to $\mathcal{\overline{W}}_{3}$ or to $\mathcal{\overline{W}}_{6}$.
\end{prop}

\subsection{An example in $\mathcal{W}_{2}$}
Let $(M, g, P)$ be a manifold with
\begin{equation}\label{prw2}
\begin{split}
    A&=\exp{(x^{1}+x^{2}-x^{3}-x^{4})}, \quad B=\sinh(x^{1}+x^{2}-x^{3}-x^{4}),\\
     C&=\exp{(x^{3}+x^{4}-x^{1}-x^{2})},
    \end{split}
\end{equation}
where $0<x^{1}+x^{2}-x^{3}-x^{4}<\ln\sqrt{3}.$

The inequalities $A>C>B>0$ are satisfied. The functions from \eqref{prw2} give a solution to \eqref{w2} but do not give a solution to \eqref{w0}.

Consequently we have the following
\begin{prop}
 The manifold $(M, g, P)$ with \eqref{prw2} belongs to $\mathcal{W}_{2}$ but does not belong to $\mathcal{W}_{0}$.
\end{prop}

\section*{Acknowledgments}
This work was partially supported by project FP17-FMI-008 of the Scientific Research Fund, Paisii Hilendarski University of Plovdiv, Bulgaria.


\end{document}